\documentstyle[11pt]{article}
\newfam\msbfam
\font\tenmsb=msbm10    \textfont\msbfam=\tenmsb
\font\sevenmsb=msbm7   \scriptfont\msbfam=\sevenmsb
\font\fivemsb=msbm5    \scriptscriptfont\msbfam=\fivemsb
\def\Bbb{\fam\msbfam \tenmsb}

\def\rr{{\Bbb R}}
\def\rz{{{\rr}^n}}

\def\fz{\infty}
\def\az{\alpha}

\def\ez{\epsilon}

\def\tz{\theta}

\def\vz{\varphi}
\def\dz{\delta}
\def\gz{\gamma}

\def\lz{\lambda}

\def\wt{\widetilde}
\def\wz{\omega}
\def\l{\left}
\def\r{\right}

\def\dsum{\displaystyle\sum}

\def\dint{\displaystyle\int}
\def\dfrac{\displaystyle\frac}
\def\dsup{\displaystyle\sup}

\def\dinf{\displaystyle\inf}

\newtheorem{thm}{\hskip\parindent Theorem}

\newtheorem{lem}{\hskip\parindent Lemma}

\newtheorem{prop}{\hskip\parindent Proposition}

\newtheorem{cor}{\hskip\parindent Corollary}

\begin{document}

\baselineskip=15pt
\renewcommand{\arraystretch}{2}
\arraycolsep=1.2pt

\title{Weighted norm inequalities for commutators of Littlewood-Paley functions related to Schr\"odinger operators
 \footnotetext{ \hspace{-0.65
cm} 2000 Mathematics Subject  Classification: 42B25, 42B20.\\
The  research was supported  by the NNSF (10971002) of China.}} %The project was supported by  the NNSF(10371004) of China.}}

\author{ Lin Tang }
\date{}
\maketitle

{\bf Abstract}\quad  Let $L=-\Delta+V$ be a Schr\"{o}dinger
operator, where $\Delta $ is the Laplacian operator on $\rz$,
while the nonnegative potential $V$ belongs to certain reverse H\"{o}lder
class. In this paper, we establish  some weighted norm inequalities
for  commutators of Littlewood-Paley functions related to Schr\"odinger operators.

\bigskip

\begin{center}{\bf 1. Introduction }\end{center}

In this paper, we consider the  Sch\"odinger differential operator
$$L=-\Delta+V(x)\ {\rm on}\ \rz,\ n\ge 3,$$where  $V$
is a  nonnegative potential satisfying  certain reverse H\"older
class.

 We say a nonnegative locally $L^q$
integral function $V(x)$ on $\rr^n$ is said to belong to
$B_q(1<q\le \fz)$ if there exists $C>0$ such that the reverse
H\"older inequality
$$\l(\dfrac 1{|B(x,r)|}\dint_{B(x,r)} V^q(y)dy\r)^{1/q}\le C\l(\dfrac
1{|B(x,r)|}\dint_{B(x,r)} V(y)dy\r) \eqno(1.1)$$ holds for every
$x\in\rz$ and $0<r<\fz$, where $B(x,r)$ denotes the ball centered
at $x$ with radius $r$. In particular, if $V$ is a nonnegative
polynomial, then $V\in B_\fz$.  It is worth pointing out that the
$B_q$ class is that, if $V\in B_q$ for some $q>1$, then there
exists $\ez>0$, which depends only $n$ and the constant $C$ in
(1.1), such that $V\in B_{q+\ez}$. Throughout this paper, we
always assume that $0\not\equiv V\in B_{n/2}$.

The study of schr\"odinger operator $L=-\triangle+V$ recently
attracted much attention; see \cite{bhs1,bhs2,dz,dz1,s1,z}.
In particular, it should be pointed out that Shen \cite{s1} proved
the Schr\"odinger type operators, such as
$\nabla(-\Delta+V)^{-1}\nabla$, $\nabla(-\Delta+V)^{-1/2}$,
$(-\Delta+V)^{-1/2}\nabla$  with $V\in B_n$, $(-\Delta+V)^{i\gz}$  with $\gz\in\rr$ and $V\in B_{n/2}$,  are standard Calder\'on-Zygmund operators.

 Recently, Bongioanni, etc,  \cite{bhs1} proved
$L^p(\rz)(1<p<\fz)$ boundedness for commutators of Riesz transforms associated with Schr\"odinger operator with $BMO(\rho)$ functions which include the class $BMO$ function, and in \cite{ bhs2} established the weighted boundedness for
Riesz transforms, fractional integrals  and Littlewood-Paley functions associated with Schr\"odinger operator with weight $A_p^{\rho}$ class which includes the Muckenhoupt weight class.  Very recently, the author \cite{t} established  the weighted norm inequalities for some Schr\"odinger type operators, which include Riesz transforms and fractional integrals  and their commutators.

In this paper, we  will continue to study weighted norm inequalities for commutators of Littlewood-Paley functions related to Schr\"odinger operators.  More precisely, we have the following results.

\begin{thm}\label{t1.1.}\hspace{-0.1cm}{\rm\bf 1.1.}\quad
Let   $1< p<\fz$. If $b\in BMO(\rho)$(defined in Section 2),
$\wz\in A_p^\rho$(defined in Section 2), then there exists a constant $C$ such that
$$\|g_b (f)\|_{L^p(\wz)}
\le C \|b\|_{BMO(\rho)}\|f\|_{L^p(\wz)}.$$
\end{thm}
where the Littlewood-Paley $g$ function related to Schr\"odinger operators is defined by
$$g(f)(x)=\l(\dint_0^\fz\l|\frac d{dt}e^{-tL}(f)(x)\r|^2tdt\r)^{1/2},\eqno(1.2)$$ and the commutator $g_b$ of $g$ with $b\in BMO(\rho)$ is defined by
$$g_b(f)(x)=\l(\dint_0^\fz\l|\frac d{dt}e^{-tL}((b(x)-b(\cdot))f)(x)\r|^2tdt\r)^{1/2}.\eqno(1.3)$$
In addition, we denote $g^*(f)(x)$ and $g^*_b(f)(x)$ in (1.2) and (1.3) if $L=\triangle$

The weighted weak-type endpoint estimate for the commutator is the
 following.

\begin{thm}\label{t1.2.}\hspace{-0.1cm}{\rm\bf 1.2.}\quad Let $b\in
BMO(\rho)$ and $\wz\in A_1^{\rho}$. There exists a constant $C>0$ such that for any $\lz>0$
$$\wz(\{x\in\rz: \ |g_bf(x)|>\lz\})\le C \dint_\rz\dfrac
{|f(x)|}\lz\l(1+\log^+\l(\frac {|f(x)|}\lz\r)\r)\wz(x)dx.$$
\end{thm}
Throughout this paper, we let $C$ denote   constants that are
independent of the main parameters involved but whose value may
differ from line to line. By $A\sim B$, we mean that there exists
a constant $C>1$ such that $1/C\le A/B\le C$.

\begin{center} {\bf 2. Preliminaries}\end{center}

We first recall some notation.  Given $B=B(x,r)$
and $\lz>0$, we will write $\lz B$ for the $\lz$-dilate ball,
which is the ball with the same center $x$ and with radius $\lz
r$. Similarly, $Q(x,r)$ denotes the cube centered at $x$ with the
sidelength $r$ (here and below only cubes
with sides parallel to the coordinate axes are considered), and $\lz Q(x,r)=Q(x,\lz r)$.
 Given a Lebesgue
measurable set $E$ and a weight $\wz$, $|E|$ will denote the
Lebesgue measure of $E$ and $\wz(E)=\int_E\wz dx$.
$\|f\|_{L^p(\wz)}$ will denote $(\int_\rz |f(y)|^p\wz(y)dy)^{1/p}$
for $0< p<\fz$.

 The function $m_V(x)$ is defined by
$$\rho(x)=\dfrac 1{m_V(x)}=\dsup_{r>0}\l\{r:\ \dfrac
{1}{r^{n-2}}\dint_{B(x,r)}V(y)dy\le 1\r\}.$$Obviously,
$0<m_V(x)<\fz$ if $V\not=0$. In particular, $m_V(x)=1$ with $V=1$
and $m_V(x)\sim (1+|x|)$ with $V=|x|^2$.

\begin{lem}\label{l2.1.}\hspace{-0.1cm}{\rm\bf 2.1(\cite{s1}).}\quad
There exists $l_0>0$ and $C_0>1$such that $$\dfrac
1{C_0}\l(1+|x-y|m_V(x)\r)^{-l_0}\le \dfrac{m_V(x)}{m_V(y)}\le
C_0\l(1+|x-y|m_V(x)\r)^{l_0/(l_0+1)}.
$$ In particular, $m_V(x)\sim m_V(y)$ if $|x-y|<C/m_V(x)$.
\end{lem}
In this paper, we write $\Psi(B)=(1+rm_V(B))^\tz$ where
$m_V(B)=\frac 1{|B|}\int_B m_V(x)dx$ and $\tz>0$, and $r$ denotes
the radius of $B$.

Obviously,
$$\Psi(B)\le \Psi(2B)\le 2^\tz \Psi(B).\eqno(2.1)$$

A weight will always mean a positive function which is locally
integrable. As \cite{bhs2}, we say that a weight $\wz$ belongs to the class
$A_p^\rho$ for $1<p<\fz$, if there is a constant $C$ such that for
all ball
  $B=B(x,r)$
$$\l(\dfrac 1{\Psi(B)|B|}\dint_B\wz(y)\,dy\r)
\l(\dfrac 1{\Psi(B)|B|}\dint_B\wz^{-\frac 1{p
-1}}(y)\,dy\r)^{p-1}\le C.$$
 We also
say that a  nonnegative function $\wz$ satisfies the $A_1^\rho$
condition if there exists a constant $C$ for all balls $B$
$$M_V(\wz)(x)\le C \wz(x), \ a.e.\ x\in\rz.$$
 where
$$M_Vf(x)=\dsup_{x\in B}\dfrac 1{\Psi(B)|B|}\dint_B|f(y)|\,dy.$$
When $V=0$, we
denote $M_{0}f(x)$ by $Mf(x)$( the standard Hardy-Littlewood
maximal function). It is easy to see that $|f(x)|\le M_{V} f(x)\le
Mf(x)$ for $a.e.\ x\in\rz$.

  We denote $A_\fz^\rho=\bigcup_{p\ge 1}A_p^\rho$. Since $\Psi(B)\ge 1$, obviously,
$A_p\subset A_p^\rho$ for $1\le p<\fz$, where $A_p^\rho$ denote
the classical Muckenhoupt weights; see \cite{gr} and \cite{m}. We
will see that $A_p\subset\subset A_p^\rho$ for $1\le p<\fz$ in
some cases.  In fact, let $\tz>0$ and $0\le\gz\le\tz$, it is easy
to check that $\wz(x)=(1+|x|)^{-(n+\gz)}\not\in A_\fz$ and
$\wz(x)dx$ is not a doubling measure, but
$\wz(x)=(1+|x|)^{-(n+\gz)}\in A_1^\rho$ provided that  $V=1$ and
$\Psi(B(x_0,r))=(1+r)^\tz$.

From the definition of $A_p^\rho$ for $1\le p<\fz$, it is easy to
see that
\begin{lem}\label{l2.2.}\hspace{-0.1cm}{\rm\bf 2.2.}\quad Let $1\le p<
\fz$. Then
\begin{enumerate}
\item[(i)]If $ 1\le p_1<p_2<\fz$, then $A_{p_1}^\rho\subset
A_{p_2}^\rho$. \item[(ii)] $\wz\in A_p^\rho$ if and only if
$\wz^{-\frac 1{p-1}}\in A_{p'}^\rho$, where $1/p+1/p'=1.$
\end{enumerate}
\end{lem}
Bongioanni, etc, \cite{bhs1} introduce a new space $BMO(\rho)$ defined by
$$\|f\|_{BMO(\rho)}=\dsup_{B\subset
\rz}\dfrac 1{\Psi(B)|B|}\dint_B|f(x)-f_B|dx<\fz,$$ where
$f_B=\frac 1{|B|}\int_Bf(y)dy$ and $\Psi(B)=(1+r/\rho(x_0))^\tz$, $B= B(x_0,r)$ and $\tz> 0$.

In particularly, Bongioanni, etc, \cite{bhs1} proved the following result for $BMO(\rho)$.
\begin{lem}\label{l2.3.}\hspace{-0.1cm}{\rm\bf 2.3.}\quad
Let $\tz>0$ and $1\le s<\fz$. If $b\in BMO(\rho)$, then
$$\l(\dfrac 1{|B|}\dint_B|b-b_B|^s\r)^{1/s}\le C_{\tz,s}\|b\|_{BMO(\rho)}\l(1+\dfrac r{\rho(x)}\r)^{\tz'},$$ for all $B=B(x,r)$, with
$x\in\rz$ and $r>0$, where $\tz'=(l_0+1)\tz$.
\end{lem}

Obviously, the classical $BMO$ is properly contained in $BMO(\rho)$; more examples see \cite{bhs1}.

From Lemma 2.3, the author \cite{t} proved the John-Nireberg inequality for $BMO(\rho)$.
\begin{prop}\label{p2.1.}\hspace{-0.1cm}{\rm\bf 2.1.}\quad
 Suppose that $f$ is in $BMO(\rho)$. There exist
positive constants $\gz$ and $C$ such that
$$\dsup_B\dfrac 1{|B|}\dint_B\exp\l\{\frac
\gz{\|f\|_{BMO(\rho)}\Psi_{\tz'}(B)}|f(x)-f_B|\r\}\,dx\le C,$$
where
$f_B=\frac 1{|B|}\int_Bf(y)dy$ and $\Psi_{\tz'}(B)=(1+r/\rho(x_0))^{\tz'}$, $B= B(x_0,r)$ and $\tz'=(l_0+1)\tz$.
\end{prop}

We remark that  balls can be replaced by cubes  in definitions of
$A_p^{\rho}$, $BMO(\rho)$ and $M_{V}$ by (2.1).

The dyadic maximal operator $M^\triangle_{V} f(x)$ is defined by
$$M_{V}^\triangle f(x):=\dsup_{x\in   Q(dyadic\ cube)}\dfrac
1{\Psi(Q)|Q|}\dint_Q|f(x)|\,dx.$$
The dyadic sharp maximal operator $M^\sharp_{V}f(x)$ is defined by
$$\begin{array}{cl}
M_{V}^\sharp f(x)&:=\dsup_{x\in Q,r<\rho(x_0)}\dfrac
1{|Q|}\dint_{Q_{x_0}}|f(y)-f_Q|\,dy+ \dsup_{x\in Q,r\ge \rho(x_0)}\dfrac
1{\Psi(Q)|Q|}\dint_{Q_{x_0}}|f|\,dy\\
&\simeq\dsup_{x\in Q,r< \rho(x_0)}\dinf_C \frac
1{|Q|}\int_{Q_{x_0}}|f(y)-C|\,dy+ \dsup_{x\in Q,r\ge \rho(x_0)}\dfrac 1{\Psi(Q) |Q|}\dint_{Q_{x_0}}|f|\,dx
\end{array}$$
where $Q_{x_0}$ denotes dyadic cubes $Q(x_0,r)$ and $f_Q=\frac
1{|Q|}\int_Q f(x)dx$.

A variant of dyadic maximal operator and dyadic sharp maximal
operator $$M^\triangle_{\dz,V}
f(x)=M^\triangle_{V}(|f|^\dz)^{1/\dz}(x)$$ and
$$M_{\dz,V}^\sharp
f(x)=M^\sharp_{V}(|f|^\dz)^{1/\dz}(x),$$ which will become the
main tool in our scheme.

In \cite{t}, the author proved the following Lemmas.
\begin{thm}\label{t2.1.}\hspace{-0.1cm}{\rm\bf 2.1.}\quad
Let $\wz\in A_\fz^\rho$. Then there  exist constant $ C, \dz^1$ such that for a locally integrable function $f$, and
for $b$ and $\gz$ positive $\gz<b<b_0=(8nC_0)^{-(l_0+2)\az}$, we
have the following inequality
$$\wz(\{x\in \rz: M^\triangle_{V} f(x)>\lz, M^\sharp_{V}f(x)\le\gz
\lz\})\le Ca^{\dz_1}\wz(\{x\in \rz:\  M^\triangle_{V}
f(x)>b\lz\})\eqno(2.1)$$ for all $\lz>0$, where $a=2^n\gz/(1-\frac
b{b_0})$.
\end{thm}

As a consequence of Theorem 2.1, we have the following result.

\begin{cor}\label{c2.1.}\hspace{-0.1cm}{\rm\bf 2.1.}\quad
Let $0<p,\  \dz<\fz$ and $\wz\in A_\fz^\rho$. There exists a
positive constant $C$ such that
$$\dint_\rz M^\triangle_{\dz,V}
 f(x)^p\wz(x)dx\le C \dint_\rz M^\sharp_{\dz,V}
f(x)^p\wz(x)dx.$$ Let  $\vz: (0,\fz)\to(0,\fz)$ be a doubling
function. Then there exists a positive constant $C$ such that
$$\sup_{\lz>0}\vz(\lz)\wz(\{x\in\rz:\ M^\triangle_{\dz,V}f(x)
>\lz\})\le C  \sup_{\lz>0}\vz(\lz)\wz(\{x\in\rz:\ M^\sharp_{\dz,V}f(x)
>\lz\})$$ for any smooth function $f$
for which the left handside is finite.
\end{cor}

\begin{prop}\label{p2.2.}\hspace{-0.1cm}{\rm\bf 2.2(\cite{t}).}\quad
Let $1<p <\fz$ and suppose that $\wz\in A_p^\rho$. If $p<p_1<\fz$,
then the equality
$$\dint_\rz|M_Vf(x)|^{p_1}\wz(x)dx\le
C_{p}\dint_\rz|f(x)|^{p_1}\wz(x)dx.$$ Further, let $1\le p<\fz$,
$\wz\in A_p^\rho$   if and only if
$$\wz(\{x\in\rz:\ M_{V}f(x)>\lz\})\le \dfrac
{C_p}{\lz^p}\dint_\rz|f(x)|^p\wz(x)dx.$$
\end{prop}
From proposition 4.1, we know that $M_V$ may be not bounded on
$L^p(\wz)$ for all $\wz\in A_p^\rho$ and $1<p<\fz$. We now need to
define a variant maximal operator $M_{V,\eta}$ for $0<\eta<\fz$ as
follows
$$M_{V,\eta}f(x)=\dsup_{x\in B}\dfrac 1{(\Psi(B))^{\eta}|B|}\dint_B|f(y)|\,dy.$$
\begin{thm}\label{t2.2.}\hspace{-0.1cm}{\rm\bf 2.2(\cite{t}).}\quad
Let $1<p<\fz$, $p'=p/(p-1)$ and suppose that $\wz\in A_p^\rho$.
 There
exists a constant $C>0$ such that
$$\|M_{V,p'}f\|_{L^p(\wz)}\le C\|f\|_{L^p(\wz)}.$$

\end{thm}

We next recall some basic definitions and facts about Orlicz
spaces, referring to \cite{r} for a complete account.

 A function $B(t): [0,\fz)\to [0,\fz)$ is called a Young function if
it is continuous, convex, increasing and satisfies $\Phi(0)=0$ and
$B\to \fz$ as $t\to\fz$. If $B$ is a Young function, we define the
$B$-average of a function $f$ over a cube $Q$ by means of the
following Luxemberg norm:
$$\|f\|_{B,Q}=\dinf\l\{\lz>0:\ \dfrac
1{|Q|}\dint_QB\l(\dfrac {|f(y)|}\lz\r)\,dy\le 1\r\}.$$
 If $A,\ B$
and $C$ are Young functions such that
$$A^{-1}(t)B^{-1}(t)\le C^{-1}(t),$$
where $ A^{-1}$ is the complementary Young function associated to
$A$, then
$$\|fg\|_{C,R}\le 2\|f\|_{A,R}\|g\|_{B,R}.$$
The examples to be considered in our study will be
$A^{-1}(t)=\log(1+t),\ B^{-1}(t)=t/\log(e+t)$ and $C^{-1}(t)=t$.
Then $A(t)\sim e^t$ and $B(t)\sim t\log(e+t)$, which gives the
generalized H\"older's inequality
$$\dfrac
1{|Q|}\dint_Q|fg|\,dy\le \|f\|_{A,Q}\|g\|_{B,Q}$$ holds. For these
example and using Theorem 2.1, if $b\in BMO(\rho)$ and $b_Q$ denotes
its average on the cube $Q$, then
$$\|(b-b_Q)/\Psi_{\tz'}(Q)\|_{expL, Q}\le C\|b\|_{BMO(\rho)}.$$
where $\tz'=(1+l_0)\tz$.

And we define the corresponding maximal function
$$M_B f(x)=\dsup_{Q:x\in Q}\|f\|_{B,Q}$$
and
$$M_{V,B} f(x)=\dsup_{Q:x\in Q}\Psi(Q)^{-1}\|f\|_{B,Q}.$$

\begin{center} {\bf 3. Some Lemmas }\end{center}
Bongioanni, etc, \cite{bhs2} proved the following result.
\begin{lem}\label{l3.1.}\hspace{-0.1cm}{\rm\bf 3.1.}\quad
Let $g^*_{loc}(f)(x)=g^*(f\chi_{B(x,\rho(x))})(x)$. Let $1<p <\fz$ and suppose that $\wz\in A_p^\rho$. Then
$$\dint_\rz|g^*_{loc}(f)(x)|^p\wz(x)dx\le
C\dint_\rz|f(x)|^p\wz(x)dx.$$ Furthermore, suppose that $\wz\in
A_1^\rho$. Then, there exists a constant $C$ such that for all
$\lz>0$
$$\wz(\{x\in\rz:\ g^*_{loc}(f)(x)>\lz\})\le \dfrac
C\lz\dint_\rz|f(x)|\wz(x)dx.$$
\end{lem}

\begin{lem}\label{l3.2.}\hspace{-0.1cm}{\rm\bf 3.2.}\quad
Let $b\in BMO(\rho)$, and $(l_0+1)\le\eta<\fz$. Set $g^*_{loc,b}(f)(x)=g^*((b(x)-b(\cdot))f\chi_{B(x,\rho(x))})(x)$. Let
$0<2\dz<\ez<1$, then
$$M^\sharp_{\dz,\eta}(g^*_{loc,b}(f))(x)\le C \|b\|_{BMO(\rho)}(M^\triangle_{\ez,\eta}(g^*_{loc}( f))(x)+M_{L\log L,V,\eta}(f)(x)), \ {\rm a.e\ }\ x\in\rz,\eqno(3.1)$$
holds for any $f\in C_0^\fz(\rz)$.
\end{lem}
{\it Proof.}\quad
We fix $x\in\rz$ and let $x\in Q=Q(x_0,r)$(dyadic cube).
To prove (3.1), we consider two cases about $r$, that is, $r<
\rho(x_0)$ and $r\ge \rho(x_0)$.

Case 1. when  $r< \rho(x_0)$. Decompose $f=f_1+f_2$, where $f_1=f\chi_{\bar Q}$,  where $
\bar Q=Q(x_0,4\sqrt{n}r)$. Let $\lz$ be a constant and $C_Q$ a constant to be
fixed along the proof. Since $0<\dz<1$,  we then have
$$\begin{array}{cl}
\l(\dfrac
1{|Q|}\dint_Q|\r.&\l.|g^*_{loc,b}(f)(y)|^\dz-|C_Q|^\dz|\,dy\r)^{1/\dz}\\&
\le \l(\dfrac
1{|Q|}\dint_Q|g^*_{loc,b}(f)(y)-C_Q|^\dz\,dy\r)^{1/\dz}\\
&\le C\l(\dfrac
1{|Q|}\dint_Q|(b(y)-\lz)g^*_{loc} f(y)|^\dz\,dy\r)^{1/\dz}\\
&\qquad+C\l(\dfrac
1{|Q|}\dint_Q|g^*_{loc}((b-\lz)f_1)(y)|^\dz\,dy\r)^{1/\dz}\\
&\qquad+C\l(\dfrac
1{|Q|}\dint_Q|g^*_{loc}((b-\lz)f_2)(y)-C_Q|^\dz\,dy\r)^{1/\dz}\\
&:=I+II+III.
\end{array}$$
 To deal with $I$, we first fix
$\lz=b_{\bar Q},$ the average of $b$ on $\bar Q$. Then for any
$1<\gz<\ez/\dz$,  note that $m_V(x)\sim m_V(x_0)$ for any $x\in \bar
Q$ and $\Psi(\bar Q)\sim 1$, by  Lemma 2.3, we then obtain
$$\begin{array}{cl}
I &\le C \l(\dfrac 1{|\bar Q|}\dint_{\bar Q}|b(y)-b_{\bar Q}|^{\dz
\gz'}\,dy\r)^{\gz'/\dz}\l(\dfrac
1{|Q|}\dint_Q|g^*_{loc}(f)(y)|^{\dz \gz}\,dy\r)^{\dz \gz}\\
&\le C \|b\|_{BMO(\rho)}M^\triangle_{\ez,\eta}(g^*_{loc}( f))(x),
\end{array}\eqno(3.2)$$ where $1/\gz'+1/\gz=1$.

For II, note that  $m_V(x)\sim m_V(x_0)$ for any $x\in \bar Q$ and $\Psi(\bar
Q)\sim 1$, by Kolmogorov's inequality and and Proposition 2.1 and Lemma 3.1, we then have
$$\begin{array}{cl}
II&\le \dfrac C{|Q|}\|g((b-b_{\bar Q}) f_1)\|_{L^{1,\fz}}\\
&\le \dfrac C{|\bar Q|}\dint_{\bar Q}|(b-b_{\bar Q})f(y)|\,dy\\
&\le
CM_{L\log L,V,\eta}f(x).
\end{array}\eqno(3.3)$$
For III, we first  fix the value of $C_Q$ by taking
$C_Q=g^*_{loc}((b-b_{\bar Q})f_2)(y_0)$  with $y_0\in Q$. Let
$b_{Q_k}=b_{Q(x_0,2^{k+1}r)}$. By Proposition 2.1, we
then obtain
$$\begin{array}{cl}
II&\le \dfrac C{|Q|}\dint_Q|g^*_{loc}((b-b_{\bar Q})f_2)(y)-g^*_{loc}((b-b_{\bar Q})f_2)(y_0)|\,dy\\
&\le \dfrac C{|Q|}\dint_Q\l[\dint_0^\fz\l(\dint_{2r<|z-x_0|\le c\rho(x_0)}|f(z)||b(z)-b_{\bar Q}|\dfrac{(t^{-n/2}|y-y_0|/\sqrt{t})}{(1+|z-y_0|/\sqrt{t})^{n+2}}dz\r)^2 tdt\r]^{1/2}dy\\
&\le \dfrac C{|Q|}\dint_Q\l[\dint_0^\fz\l(\dint_{2r<|z-x_0|\le c\rho(x_0)}|f(z)||b(z)-b_{\bar Q}|\dfrac{rt}{(t+|z-x_0|)^{n+2}}dz\r)^{2}
\dfrac{dt} t\r]^{1/2}dy\\
&\le \dfrac C{|Q|}\dint_Q\l[\dint_{2r<|z-x_0|\le c\rho(x_0)}r|f(z)||b(z)-b_{\bar Q}|\l(\dint_0^\fz\dfrac{t}{(t+|z-x_0|)^{2(n+2)}}dt\r)^{1/2}
dz\r]dy\\
&\le \dfrac C{|Q|}\dint_Q\l[\dint_{2r<|z-x_0|\le c\rho(x_0)}r|f(z)||b(z)-b_{\bar Q}||z-x_0|^{-(n+1)}
dz\r]dy\\
&\le \dfrac C{|Q|}\dint_Q\l[\dsum_{k=1}^{k_0}\dfrac{2^{-k}}{(2^kr)^n}\dint_{|z-x_0|\le 2^{k+1}}|f(z)||b(z)-b_{\bar Q}|
dz\r]dy\\
&\le C\|b\|_{BMO(\rho)}M_{L\log L,V,\eta}(f)(x),
\end{array}\eqno(3.4)$$
where the integer $k_0$ satisfies  $2^{k_0}r\le c\rho(x_0)\le 2^{k_0+1}$ and $c=C_0n2^{l_0+4}$.

Case 2. When  $r\ge \rho(x_0)$.  Decompose $f=f_1+f_2$, where $f_1=f\chi_{\bar Q}$,  where $
\bar Q=Q(x_0,C_02^{l_0+4}\sqrt{n}r)$. Since $0<2\dz< \ez<1$, so
 $a=\eta/\dz$ and $\ez/\dz>2$, then
$$\begin{array}{cl}
\dfrac 1{\Psi(Q)^a}\l(\dfrac
1{|Q|}\dint_Q\r.&\l.|g^*_{loc,b}(f))f(y)|^\dz\,dy\r)^{1/\dz}\\
 & \le
\dfrac 1{\Psi(Q)^a}\l(\dfrac
1{|Q|}\dint_Q|(b(y)-\lz)g^*_{loc}(f)(y)+g^*_{loc}((b-\lz)f)(y)|^\dz\,dy\r)^{1/\dz}\\
&\le C\dfrac 1{\Psi(Q)^a}\l(\dfrac
1{|Q|}\dint_Q|(b(y)-\lz)g^*_{loc}(f))(y)|^\dz\,dy\r)^{1/\dz}\\
&\qquad+C\dfrac 1{\Psi(Q)^a}\l(\dfrac
1{|Q|}\dint_Q|g^*_{loc}((b-\lz)f_1)(y)|^\dz\,dy\r)^{1/\dz}\\
&\qquad+C\dfrac 1{\Psi(Q)^a}\l(\dfrac
1{|Q|}\dint_Q|g^*_{loc}((b-\lz)f_2)(y)|^\dz\,dy\r)^{1/\dz}\\
&:=I+II+III.
\end{array}$$
 To deal with $I$, we first fix
$\lz=b_{\bar Q},$ the average of $b$ on $\bar Q$. Then for any
$2\le\gz<\ez/\dz$, note that $l_0+1\le\eta$, by Lemma  2.3,  we then have
$$\begin{array}{cl}
I &\le C\dfrac 1{\Psi_{\tz'}(Q)} \l(\dfrac 1{|\bar Q|}\dint_{\bar
Q}|b(y)-b_{\bar Q}|^{\dz \gz'}\,dy\r)^{1/(r'\dz)}\\
&\qquad\times\dfrac
{\Psi_{\tz'}(Q)} {\Psi(Q)^{a-\eta/(2\dz)}}\l(\dfrac
1{\Psi(Q)^\eta|Q|}\dint_Q|g^*_{loc}(f))(y)|^{\dz \gz}\,dy\r)^{1/(\dz \gz)}\\
\end{array}$$
$$\begin{array}{cl}
&\le C \|b\|_{BMO(\rho)}M^\triangle_{\ez,\eta}(g^*_{loc}(f)) f)(x),
\end{array}\eqno(3.5)$$ where $1/\gz'+1/\gz=1$.

For II, we recall that $g^*_{loc}$ is weak type $(1,1)$ by Lemma 3.1. By
Kolmogorov's inequality and Proposition 2.1, we then have

$$\begin{array}{cl}
II&\le \dfrac C{\Psi(Q)^a} \dfrac 1{|Q|}\|g^*_{loc}((b-b_{\bar Q}) f_1)\|_{L^{1,\fz}}\\
&\le \dfrac C{\Psi(Q)^{a}}\dfrac 1{|\bar Q|}\dint_{\bar Q}|(b-b_{\bar Q})f(y)|\,dy\\
&\le CM_{L\log L,V,\eta}f(x).
\end{array}\eqno(3.6)$$
Finally, for III, notice that $B(y,\rho(y))\subset Q(x_0,C_02^{l_0+4}\sqrt{n}r)$ for any $y\in Q$, then $III=0$.

 From  (3.2)--(3.6), we get (3.1). Hence the proof is
 finished.   \hfill$\Box$

We next consider several maximal operators, which play an important role in this paper.
$$M_{V,\eta}f(x)=\dsup_{x\in B}\dfrac 1{(\Psi(B))^{\eta}|B|}\dint_B|f(y)|\,dy,$$
$$\wt M^b_{V,\eta}f(x)=\dsup_{\ez>0}\dfrac 1{(1+ \ez\psi(B(x,\ez)))^{\tz\eta}}\dint_\rz\ez^{-n}\vz(\frac{x-y}\ez)|f(y)|dy,$$ and their commutators
$$M^b_{V,\eta}f(x)=\dsup_{x\in B}\dfrac 1{(\Psi(B))^{\eta}|B|}\dint_B|b(x)-b(y)||f(y)|\,dy,$$
$$\wt M^b_{V,\eta}f(x)=\dsup_{\ez>0}\dfrac 1{(1+ \ez\psi(B(x,\ez)))^{\tz\eta}}\dint_\rz\ez^{-n}\vz(\frac{x-y}\ez)|b(x)-b(y)||f(y)|dy,$$
where $\psi(B(x,\ez))=\frac 1{B(x,\ez)}\int_{B(x,\ez)}\rho(y)^{-1}dy$

Obviously, we have $$M^b_{V,\eta'}f(x)\le C \wt M^b_{V,\eta}f(x),\eqno(3.7)$$ where $\eta'=(l_0+1)\eta$ and $ \eta>0$.
\begin{lem}\label{l3.3.}\hspace{-0.1cm}{\rm\bf 3.3.}\quad
Let $b\in BMO(\rho)$, and $(l_0+1)(1+1/\tz)\le\eta<\fz$, $\eta_1=(l_0+1)\eta$ and $\eta_2=(l_0+1)\eta_1(1+1/\tz)$.  Let
$0<2\dz<\ez<1$, then
$$M^\sharp_{\dz,\eta}(\wt M^b_{V,\eta_2}(f))(x)\le C \|b\|_{BMO(\rho)}(M^\triangle_{\ez,\eta}(\wt M_{V,\eta_2}( f))(x)+M_{L\log L,V,\eta}(f)(x)), \ {\rm a.e\ }\ x\in\rz,\eqno(3.9)$$
holds for any $f\in C_0^\fz(\rz)$.
\end{lem}
{\it Proof.}\quad
We fix $x\in\rz$ and let $x\in Q=Q(x_0,r)$(dyadic cube).
To prove (3.9), we consider two cases about $r$, that is, $r<
\rho(x_0)$ and $r\ge \rho(x_0)$.

Case 1. when  $r< \rho(x_0)$. Decompose $f=f_1+f_2$, where $f_1=f\chi_{\bar Q}$,  where $
\bar Q=Q(x_0,4\sqrt{n}r)$. Let $\lz$ be a constant and $C_Q$ a constant to be
fixed along the proof. Since $0<\dz<1$,  we then have

$$\begin{array}{cl}
\l(\dfrac
1{|Q|}\dint_Q|\r.&\l.|\wt M^b_{V,\eta_2}(f)(y)|^\dz-|C_Q|^\dz|\,dy\r)^{1/\dz}\\&
\le \l(\dfrac
1{|Q|}\dint_Q|\wt M^b_{V,\eta_2}(f)(y)-C_Q|^\dz\,dy\r)^{1/\dz}\\
&\le C\l(\dfrac
1{|Q|}\dint_Q|(b(y)-\lz)\wt M_{V,\eta_2}(f)(y)|^\dz\,dy\r)^{1/\dz}\\
&\qquad+C\l(\dfrac
1{|Q|}\dint_Q|\wt M_{V,\eta_2}((b-\lz)f_1)(y)|^\dz\,dy\r)^{1/\dz}\\
&\qquad+C\l(\dfrac
1{|Q|}\dint_Q|\wt M_{V,\eta_2}((b-\lz)f_2)(y)-C_Q|^\dz\,dy\r)^{1/\dz}\\
&:=I+II+III.
\end{array}$$
 To deal with $I$, we first fix
$\lz=b_{\bar Q},$ the average of $b$ on $\bar Q$. Then for any
$1<\gz<\ez/\dz$,  note that $m_V(x)\sim m_V(x_0)$ for any $x\in \bar
Q$ and $\Psi(\bar Q)\sim 1$, by  Proposition 2.1, we then obtain
$$\begin{array}{cl}
I &\le C \l(\dfrac 1{|\bar Q|}\dint_{\bar Q}|b(y)-b_{\bar Q}|^{\dz
\gz'}\,dy\r)^{\gz'/\dz}\l(\dfrac
1{|Q|}\dint_Q|\wt M_{V,\eta_2}(f)(y)|^{\dz \gz}\,dy\r)^{\dz \gz}\\
&\le C \|b\|_{BMO(\rho)}M^\triangle_{\ez,\eta}(\wt M_{V,\eta_2}( f))(x),
\end{array}\eqno(3.10)$$ where $1/\gz'+1/\gz=1$.

For II, note that  $m_V(x)\sim m_V(x_0)$ for any $x\in \bar Q$ and $\Psi(\bar
Q)\sim 1$, by Kolmogorov's inequality and Theorem 2.1, by the weak (1,1) of $\wt M_{V,\eta_2}$, we then have
$$\begin{array}{cl}
II&\le \dfrac C{|Q|}\|\wt M_{V,\eta_2}((b-b_{\bar Q}) f_1)\|_{L^{1,\fz}}\\
&\le \dfrac C{|\bar Q|}\dint_{\bar Q}|(b-b_{\bar Q})f(y)|\,dy\\
&\le
C\|b\|_{BMO(\rho)}M_{L\log L,V,\eta}f(x).
\end{array}\eqno(3.11)$$
For III, we   fix the value of $C_Q$ by taking
$C_Q=\wt M_{V,\eta_2}((b-b_{\bar Q}) f_2))(y_0)$ for some $y_0\in Q$.
 We now estimate
$E:=|\wt M_{V,\eta}((b-b_{\bar Q}) f_2)(y)-C_Q|$ for any $y\in Q$.
$$\begin{array}{cl}
E&= \l|\dsup_{\ez>0}\dfrac 1{(1+\ez\psi(B(y,\ez)))^{\tz\eta_2}}\dint_\rz\ez^{-n}\vz(\frac{y-z}\ez)|b(z)-b_{\bar Q}||f_2(z)|dz\r.\\
&\quad - \l.\dsup_{\ez>0}\dfrac 1{(1+ \ez\psi(B(y_0,\ez)))^{\tz\eta_2}}\dint_\rz\ez^{-n}\vz(\frac{y_0-z}\ez)|b(z)-b_{\bar Q}||f_2(z)|dz\r|\\
&\le \dsup_{\ez>0}\dfrac 1{(1+\ez\psi(B(y,\ez)))^{\tz\eta_2}}\dint_\rz\ez^{-n}|\vz(\frac{y-z}\ez)-\vz(\frac{y_0-z}\ez)||b(z)-b_{\bar Q}||f_2(z)|dz\\
&\quad+\dsup_{\ez>0}\l|\dfrac 1{(1+ \ez\psi(B(y_0,\ez)))^{\tz\eta_2}}-\dfrac 1{(1+ \ez\psi(B(z,\ez)))^{\tz\eta_2}}\r|\\
&\qquad\qquad\times\dint_\rz\ez^{-n}\vz(\frac{y_0-z}\ez)|b(z)-b_{\bar Q}||f_2(z)|dz\\
\end{array}$$
$$\begin{array}{cl}
&= \dsup_{\ez>r}\dfrac 1{(1+\ez\psi(B(y,\ez)))^{\tz\eta_2}}\dint_\rz\ez^{-n}|\vz(\frac{y-z}\ez)-\vz(\frac{y_0-z}\ez)||b(z)-b_{\bar Q}||f_2(z)|dz\\
&+\dsup_{\ez>r}\l|\dfrac 1{(1+ \ez\psi(B(y_0,\ez)))^{\tz\eta_2}}-\dfrac 1{(1+ \ez\psi(B(y,\ez)))^{\tz\eta_2}}\r|\dint_\rz\ez^{-n}\vz(\frac{y_0-y}\ez)|b(z)-b_{\bar Q}||f_2(z)|dz\\
&\le \dsup_{\ez>r}\dfrac C{(1+\frac \ez{\rho(y)})^{\tz\eta_1}}\dint_{r\le|z-y|\le 8\ez}\ez^{-n}\frac r\ez|b(y)-b_{\bar Q}||f(y)|dy\\
&\quad+ C\dsup_{\ez>r}\ez|\psi(B(y_0,\ez))-\psi(B(y,\ez))|\l|\dfrac 1{(1+ \ez\psi(B(y_0,\ez)))^{\tz\eta_2}}+\dfrac 1{(1+ \ez\psi(B(y,\ez)))^{\tz\eta_2}}\r|\\
&\qquad\qquad\qquad\times\dint_\rz\ez^{-n}\vz(\frac{y_0-z}\ez)|b(z)-b_{\bar Q}||f_2(z)|dz\\
&\le \dsup_{\ez>r}\dfrac C{(1+\frac \ez{\rho(y)})^{\tz\eta_1}}\dint_{r\le|z-y|\le 8\ez}\ez^{-n}\frac r\ez|b(z)-b_{\bar Q}||f(z)|dz\\
&\quad+ C\dsup_{\ez>r}(\ez\rho(y)^{-1})^{l_0+1}\frac r\ez\dfrac 1{(1+\ez/\rho(y))^{\tz\eta_1(1+\frac 1\tz)}}\dint_\rz\ez^{-n}\vz(\frac{y_0-z}\ez)|b(z)-b_{\bar Q}||f_2(z)|dz\\
&\le \dsup_{\ez>r}\dfrac C{(1+\frac \ez{\rho(y)})^{\tz\eta_1}}\dint_{r\le|z-y|\le 8\ez}\ez^{-n}\frac r\ez|b(z)-b_{\bar Q}||f(z)|dz\\
&\quad+ C\dsup_{\ez>r}\dfrac 1{(1+\ez/\rho(y))^{\tz\eta_1}}\dint_{r\le|z-y|\le 8\ez}\ez^{-n}\frac r\ez|b(z)-b_{\bar Q}||f(z)|dz\\
&\le \dsup_{\ez>r}\dsum_{k=1}^{[\ln (\frac{8\ez}r)]+1}\dfrac C{(1+\frac \ez{\rho(y)})^{\tz\eta_1}}\ez^{-n-1}r\dint_{|z-y|\le 2^kr}|b(z)-b_{\bar Q}||f(z)|dz\\
&\le \dsup_{\ez>r}\dsum_{k=1}^{[\ln (\frac{8\ez}r)]+1}\frac r\ez\dfrac C{(1+\frac {2^kr}{\rho(y)})^{\tz\eta_1}(2^kr)^{n}}\dint_{|z-y|\le 2^kr}|b(z)-b_{\bar Q}||f(z)|dz\\
&\le \dsup_{\ez>r}\dsum_{k=1}^{[\ln (\frac{8\ez}r)]+1}\frac r\ez
k\|b\|_{BMO(\rho)}M_{L\log L,V,\eta}f(x)\\
&\le C\|b\|_{BMO(\rho)}M_{L\log L,V,\eta}f(x).
\end{array}$$
Hence,
$$III\le  C\|b\|_{BMO(\rho)}M_{L\log L,V,\eta}f(x).\eqno(3.12)$$
Case 2. when  $r> \rho(x_0)$.  Let $f_1, f_2$ be above.
  We then have
$$\begin{array}{cl}
\l(\dfrac
1{|Q|}\dint_Q|\wt M^b_{V,\eta_2}(f)(y)|^\dz\,dy\r)^{1/\dz} &\le C\l(\dfrac
1{|Q|}\dint_Q|(b(y)-\lz)\wt M_{V,\eta_2}(f)(y)|^\dz\,dy\r)^{1/\dz}\\
&\qquad+C\l(\dfrac
1{|Q|}\dint_Q|\wt M_{V,\eta_2}((b-\lz)f_1)(y)|^\dz\,dy\r)^{1/\dz}\\
&\qquad+C\l(\dfrac
1{|Q|}\dint_Q|\wt M_{V,\eta_2}((b-\lz)f_2)(y)|^\dz\,dy\r)^{1/\dz}\\
&:=I_1+II_1+III_1.
\end{array}$$
To deal with $I_1$, we first fix
$\lz=b_{\bar Q},$ the average of $b$ on $\bar Q$. Then for any
$2\le\gz<\ez/\dz$,  by  Lemma  2.3, we then obtain  that
$$\begin{array}{cl}
I &\le C\dfrac 1{\Psi_{\tz'}(Q)} \l(\dfrac 1{|\bar Q|}\dint_{\bar
Q}|b(y)-b_{\bar Q}|^{\dz \gz'}\,dy\r)^{1/(r'\dz)}\\
&\qquad\times\dfrac
{\Psi_{\tz'}(Q)} {\Psi(Q)^{a-\eta/(2\dz)}}\l(\dfrac
1{\Psi(Q)^\eta|Q|}\dint_Q|g^*_{loc}(f))(y)|^{\dz \gz}\,dy\r)^{1/(\dz \gz)}\\
&\le C \|b\|_{BMO(\rho)}M^\triangle_{\ez,\eta}(\wt M_{V,\eta_2}( f))(x),
\end{array}\eqno(3.13)$$ where $1/\gz'+1/\gz=1$.

For $II_1$,  by Kolmogorov's inequality and Proposition 2.1, by the weak (1,1) of $\wt M_{V,\eta_2}$, we then have
$$\begin{array}{cl}
II_1&\le \dfrac C{|Q|}\|\wt M_{V,\eta_2}((b-b_{\bar Q}) f_1)\|_{L^{1,\fz}}\\
&\le \dfrac C{|\bar Q|}\dint_{\bar Q}|(b-b_{\bar Q})f(y)|\,dy\\
&\le
C\|b\|_{BMO(\rho)}M_{L\log L,V,\eta}f(x).
\end{array}\eqno(3.14)$$
For $III_1,$  we have for any $y\in Q$,
$$\begin{array}{cl}
\wt M_{V,\eta_2}((b-b_{\bar Q}) f_2)(y)&= \dsup_{\ez>0}\dfrac 1{(1+\ez\psi(B(y,\ez)))^{\tz\eta_2}}\dint_\rz\ez^{-n}\vz(\frac{y-z}\ez)|b(z)-b_{\bar Q}||f_2(z)|dz\\
&= \dsup_{\ez>r}\dfrac 1{(1+\ez\psi(B(y,\ez)))^{\tz\eta_2}}\dint_\rz\ez^{-n}\vz(\frac{y-z}\ez)|b(z)-b_{\bar Q}||f_2(z)|dz\\
&\le \dsup_{\ez>r}\dfrac C{(1+\frac \ez{\rho(y)})^{\tz\eta_1}}\dint_{r\le|z-x|\le 8\ez}\ez^{-n}|b(y)-b_{\bar Q}||f(y)|dy\\
&\le \dsup_{\ez>r}\dfrac C{(1+\frac \ez{\rho(y)})^{\tz\eta_1-1}}\frac r\ez\dint_{r\le|z-x|\le 8\ez}\ez^{-n}|b(y)-b_{\bar Q}||f(y)|dy\\
&\le \dsup_{\ez>r}\dfrac C{(1+\frac \ez{\rho(y)})^{\tz(l_0+1)\eta}}\frac r\ez\dint_{r\le|z-x|\le 8\ez}\ez^{-n}|b(y)-b_{\bar Q}||f(y)|dy\\
&\le
C\|b\|_{BMO(\rho)}M_{L\log L,V,\eta}f(x).
\end{array}\eqno(3.15)$$
 From  (3.10)--(3.15), we get (3.9). Hence the proof is
 finished.      \hfill$\Box$

\begin{lem}\label{l3.4.}\hspace{-0.1cm}{\rm\bf 3.4.}\quad
Let  $2\le\eta<\fz$,  $\wz\in A_1^\rho$ and
$B(t)=t\log(e+t)$. Then there exists a constant $C>0$ such that
for all $t>0$
$$\wz(\{x\in\rz:\ M_{B,V,\eta}f(x)>t\})\le
C\dint_\rz B\l(\dfrac {|f(x)|}t\r)\wz(x)dx.\eqno(3.16)$$
\end{lem}
{\it Proof.}\quad Let $K$ be any compact subset in $\{x\in\rz: M_{L\log
L,\vz,\eta}(f)(x)>\lz\})$. For any $x\in K$, by a standard
covering lemma, it is possible to choose  cubes $Q_1,\cdots, Q_m$
with pairwise disjoint interiors such that $K\subset
\bigcup_{j=1}^m 3Q_j$ and with $\|f\|_{L\log L,\vz,Q_j}>\lz$,
$j=1,\cdots, m$. This implies
$$\Psi(Q_j)^2|Q_j|\le \dint_{Q_j}\dfrac {|f(y)|}\lz\l(1+\log^+\l(\dfrac
{|f(y)|}\lz\r)\r)\,dy.$$ From this, by (vi) in  Lemma 2.1 with
$p=1$ and $E=Q$, we obtain that
$$\begin{array}{cl}
\wz(3Q_j)&\le C\Psi(Q_j)\wz(Q_j)\\
&=C\Psi(Q_j)^2|Q_j|
\dfrac{\wz(Q_j)}{\Psi(Q_j))|Q_j|} \\
&\le C\dfrac{\wz(Q_j)}{\Psi(Q_j)|Q_j|}\dint_{Q_j}\dfrac
{|f(y)|}\lz\l(1+\log^+\l(\dfrac {|f(y)|}\lz\r)\r)\,dy\\
&\le C\dinf_{Q_j}\wz(x) \dint_{Q_j}\dfrac
{|f(y)|}\lz\l(1+\log^+\l(\dfrac {|f(y)|}\lz\r)\r)\,dy\\
&\le C\dint_{Q_j}\dfrac {|f(y)|}\lz\l(1+\log^+\l(\dfrac
{|f(y)|}\lz\r)\r)\wz(y)\,dy.
\end{array}$$
Thus, (3.16) holds, hence, the proof is complete.       \hfill$\Box$

Finally, the author \cite{t} proved the following result.
\begin{lem}\label{l3.5.}\hspace{-0.1cm}{\rm\bf 3.5.}\quad
Let $ 0<\eta<\fz$ and $M_{V,\eta/2}f$ be locally
integral. Then there exist positive constants $C_1$ and $C_2$
independent of $f$ and $x$ such that
$$C_1M_{V,\eta}M_{V,\eta+1}f(x)\le M_{L\log L, V,\eta+1}f(x)\le
C_2M_{V,\eta/2} M_{V,\eta/2}f(x).$$
\end{lem}

\begin{center} {\bf 4.  Proof of some theorems}\end{center}
{\it Proof of Theorem 1.1.}\quad  We adapt a similar argument of Theorem 5 in \cite{bhs2}. As before, we define
$$g_{loc,b}(f)(x)=g((b(x)-b(\cdot))f\chi_{B(x,\rho(x))})(x),\ \  g_{glob,b}(f)(x)=g((b(x)-b(\cdot))f\chi_{B^c(x,\rho(x))})(x).$$
Thus
$$\|g_b(f)\|_{L^p(\wz)}\le \|g_{loc,b}(f)\|_{L^p(\wz)}+\|g_{glob,b}(f)\|_{L^p(\wz)}.$$
We start with $g_{glob,b}$. Denoting by $q_t$ the kernel of $\frac{d}{dt}e^{-tL}$, from (2.7) of \cite{dz1}, for any  $N>0$, we have
$$|q_t(x,y)|\le\dfrac {C_N}{t^{n/2+1}}\l(1+\frac t{\rho(x)^2}+\frac t{\rho(y)^2}\r)^{-N}e^{-\frac{|x-y|^2}{ct}}.\eqno(4.1)$$
Hence,
$$\begin{array}{cl}
&\l|\dint_{|x-y|>\rho(x)}q_t(x,y)(b(x)-b(y))f(y)dy\r|\\
&\qquad\le Ct^{-n/2-1}\l(1+\frac t{\rho(x)^2}\r)^{-N}
\dint_{|x-y|>\rho(x)}e^{-\frac{|x-y|^2}{ct}}
|b(x)-b(y)||f(y)|dy\\
&\qquad\le Ct^{\frac {M-d}2-1}\l(1+\frac t{\rho(x)^2}\r)^{-N}
\dint_{|x-y|>\rho(x)}\dfrac{|b(x)-b(y)||f(y)|}{|x-y|^{M}}dy\\
\end{array}$$
$$\begin{array}{cl}
&\qquad\le C\dfrac{t^{\frac {M-n}2-1}}{\rho(x)^{M-n}}\l(1+\frac t{\rho(x)^2}\r)^{-N}
\dsum_{k=1}^\fz \dfrac{2^{-k(M-n-\tz\eta)}}{2^{k\tz\eta}|2^k\rho(x)|^{n}}\dint_{|x-y|<2^k\rho(x)}|b(x)-b(y)||f(y)|dy\\
&\qquad\le C\dfrac{t^{\frac {M-n}2-1}}{\rho(x)^{M-n}}M_{V,\eta}^bf(x).
\end{array}$$
Then,
$$g_{glob,b}(f)(x)\le CM_{V,\eta}^bf(x)\l(\dint_0^\fz\l(\dfrac t{\rho(x)^2}\r)^{M-n}\l(1+\dfrac t{\rho(x)^2}\r)^{-2N}\dfrac{dt}t\r)^{1/2}\le CM_{V,\eta}^bf(x).$$
Choose $M$ and $N$ such that $M-n>\tz\eta$ and $2N>M-n$. Therefore, the estimates for $g_{glob,b}$ follow from those for $M_{V,\eta}^bf(x)$ by Lemmas 3.3 and 3.5.

To deal with $g_{loc,b}$ we write
$$g_{loc,b}(f)(x)\le I(x)+g_{loc,b}^*(f)(x)+II(x),\eqno(4.2)$$ where $g_{loc,b}^*(f)(x)$ is defined in Lemma 3.2,
$$I(x)=\l(\dint_0^{\rho(x)^2}\l|\dint_{|x-y|<\rho(x)}[q_t(x,y)-\wt q_t(x,y)](b(x)-b(y))f(y)dy\r|^2tdt\r)^{1/2},$$
where $\wt q_t$ is the kernel of $\frac d{dt}e^{t\triangle}$, and
$$II(x)=\l(\dint_{\rho(x)^2}^\fz\l|\dint_{|x-y|<\rho(x)}q_t(x,y)(b(x)-b(y))f(y)dy\r|^2tdt\r)^{1/2}.$$
For $II(x)$, by (4.1) with $N=1/2$,
$$\begin{array}{cl}
II(x)&\le C\l(\dint_{\rho(x)^2}^\fz\l(\dfrac{\rho(x)}t\r)^2
\l|\dint_{|x-y|<\rho(x)}t^{-n/2}e^{-\frac{|x-y|^2}{ct}}|b(x)-b(y)||f(y)|dy\r|^2tdt\r)^{1/2}\\
&\le C\l(\dint_{\rho(x)^2}^\fz\l(\dfrac{\rho(x)}t\r)^2
\l|\rho(x)^{-n}\dint_{|x-y|<\rho(x)}|b(x)-b(y)||f(y)|dy\r|^2tdt\r)^{1/2}\\
&\le CM_{V,\eta}^bf(x)\l(\dint_{\rho(x)^2}^\fz\l(\dfrac{\rho(x)}t\r)^2dt\r)^{1/2}\\
&\le CM_{V,\eta}^bf(x).
\end{array}\eqno(4.3)$$
For $I(x)$, adapting the same argument of pages 578-579 in\cite {bhs2}, we obtain for some $\dz>0$ and $\ez>0$
$$\begin{array}{cl}
I(x)&\le C\l(\dint_0^{\rho(x)^2}
\l(\dfrac{\sqrt{t}}{\rho(x)}\r)^\dz\l|\dint_{|x-y|<\rho(x)}t^{-n/2}e^{-\ez\frac{|x-y|^2}{t}}
|b(x)-b(y)||f(y)|dy\r|^2\dfrac{dt}t\r)^{1/2}\\
&\le C\l(\dint_0^{\rho(x)^2}
\l(\dfrac{\sqrt{t}}{\rho(x)}\r)^\dz\l|\dint_{\sqrt{t}\le |x-y|<\rho(x)}t^{-n/2}e^{-\ez\frac{|x-y|^2}{t}}
|b(x)-b(y)||f(y)|dy\r|^2\dfrac{dt}t\r)^{1/2}\\
&\quad+C\l(\dint_0^{\rho(x)^2}
\l(\dfrac{\sqrt{t}}{\rho(x)}\r)^\dz\l|\dint_{|x-y|<\sqrt{t}}t^{-n/2}e^{-\ez\frac{|x-y|^2}{t}}
|b(x)-b(y)||f(y)|dy\r|^2\dfrac{dt}t\r)^{1/2}\\
\end{array}$$
$$\begin{array}{cl}
&\le C\l(\dint_0^{\rho(x)^2}
\l(\dfrac{\sqrt{t}}{\rho(x)}\r)^\dz\l|\dsum_{k=1}^{[\log_2(\rho(x)/\sqrt{t})]+1}\dint_{\sqrt{t}\le |x-y|<\rho(x)}
\dfrac{|b(x)-b(y)||f(y)|}{|x-y|^n}dy\r|^2\dfrac{dt}t\r)^{1/2}\\
&\quad+C\l(\dint_0^{\rho(x)^2}
\l(\dfrac{\sqrt{t}}{\rho(x)}\r)^\dz\l|\dint_{|x-y|<\sqrt{t}}t^{-n/2}
|b(x)-b(y)||f(y)|dy\r|^2\dfrac{dt}t\r)^{1/2}\\
&\le C\l(\dint_0^{\rho(x)^2}
\l(\dfrac{\sqrt{t}}{\rho(x)}\r)^\dz  ([\log_2(\rho(x)/\sqrt{t})]+1 )^2\dfrac{dt}t\r)^{1/2} M_{V,\eta}^bf(x)\\
&\quad+C\l(\dint_0^{\rho(x)^2}
\l(\dfrac{\sqrt{t}}{\rho(x)}\r)^\dz \dfrac{dt}t\r)^{1/2} M_{V,\eta}^bf(x)\\
&\le CM_{V,\eta}^bf(x).
\end{array}\eqno(4.4)$$
 From (4.2), (4.3) and (4.4), we can obtain the desired result by Lemmas 3.2, 3.3, 3.5 and Theorem 2.2.   \hfill$\Box$

 \bigskip

 {\it Proof Theorem 1.2.}\quad By (4.1)-(4.4) and using Lemmas 3.2, 3.3, 3.4, 3.5 and Proposition 2.2,
   by adapting an argument in \cite{p},  we can obtain the desired result.       \hfill$\Box$

Finally, we consider the maximal operator of the diffusion semi-group
$$T^*f(x)=\dsup_{t>0}e^{-tL}f(x)=\dsup_{t>0}\dint_\rz k_t(x,y)f(y)dy,$$
 and it's commutator
$$T_b^*f(x)=\dsup_{t>0}e^{-tL}f(x)=\dsup_{t>0}\dint_\rz k_t(x,y)(b(x)-b(y))f(y)dy,$$
where $k_t$ is the kernel of the operator $e^{-tL},\ t>0$.
\begin{thm}\label{t4.1.}\hspace{-0.1cm}{\rm\bf 4.1.}\quad Let $b\in BMO(\rho)$ and $T_b^*f$ be as above.
\begin{enumerate}
\item[(i)] If   $1< p<\fz$,
$\wz\in A_p^\rho$, then there exists a constant $C$ such that
$$\|T^*_b f\|_{L^p(\wz)}
\le C \|b\|_{BMO(\rho)}\|f\|_{L^p(\wz)}.$$
\item[(ii)]If $\wz\in A_1^{\rho}$, then there exists a constant $C>0$ such that for any $\lz>0$
$$\wz(\{x\in\rz: \ |T^*_bf(x)|>\lz\})\le C \dint_\rz\dfrac
{|f(x)|}\lz\l(1+\log^+\l(\frac {|f(x)|}\lz\r)\r)\wz(x)dx.$$
\end{enumerate}
\end{thm}
{\it Proof}.\quad We first recall the kernel $k_t$ has the following property (see \cite{dz1})
$$0\le k_t(x,y)\le C_Nt^{-n/2}e^{-\frac{|x-y|^2}{5t}}\l(1+\frac t{\rho(x)^2}+\frac t{\rho(y)^2}\r)^{-N}.\eqno(4.5)$$
Then
$$\begin{array}{cl}
|T^*_bf(x)|&\le \dsup_{t>0}\dint_\rz k_t(x,y)|(b(x)-b(y))f(y)|dy\\
&\le \dsup_{t>0}\dint_{|x-y|<\rho(x)} k_t(x,y)|(b(x)-b(y))f(y)|dy\\
&\quad+ \dsup_{t>0}\dint_{|x-y|\ge\rho(x)} k_t(x,y)|(b(x)-b(y))f(y)|dy\\
&:=I(x)+II(x).
\end{array}$$
For $I(x)$, by (4.5), we then have
$$\begin{array}{cl}
I(x)&\le \dsup_{0<\sqrt{t}<\rho(x)}\dint_{|x-y|<\sqrt{t}} k_t(x,y)|(b(x)-b(y))f(y)|dy\\
&\quad+\dsup_{0<\sqrt{t}<\rho(x)}\dint_{\sqrt{t}\le |x-y|<\rho(x)} k_t(x,y)|(b(x)-b(y))f(y)|dy\\
&\quad+\dsup_{\sqrt{t}\ge\rho(x)}\dint_{ |x-y|\le\rho(x)} k_t(x,y)|(b(x)-b(y))f(y)|dy\\
&\le C\dsup_{0<\sqrt{t}<\rho(x)}\dint_{|x-y|<\sqrt{t}} t^{-n/2}|(b(x)-b(y))f(y)|dy\\
&\quad+C\dsup_{0<\sqrt{t}<\rho(x)}\dint_{\sqrt{t}\le |x-y|<\rho(x)} \sqrt{t}|x-y|^{-(n+1)}|(b(x)-b(y))f(y)|dy\\
&\quad+\dsup_{\sqrt{t}\ge\rho(x)}\rho(x)^{-n}\dint_{ |x-y|<\rho(x)} |(b(x)-b(y))f(y)|dy\\
&\le CM_{V,\eta}^bf(x).
\end{array}\eqno(4. 6)$$
For $II(x)$,  by (4.5) again, we then obtain  that
$$\begin{array}{cl}
II(x)&\le \dsup_{0<t}t^{-n/2}\l(1+\frac {\sqrt{t}}{\rho(x)}\r)^{-N}\dint_{|x-y|\ge\rho(x)} e^{-\frac{|x-y|^2}{5t}}|(b(x)-b(y))f(y)|dy\\
&\le \dsup_{0<t}\l(\frac {\sqrt{t}}{\rho(x)}\r)^{M-n}\l(1+\frac {\sqrt{t}}{\rho(x)}\r)^{-N}\\
&\qquad\qquad\times\dsum_{k=1}^\fz \dfrac{2^{-k(M-n-\tz\eta)}}{2^{k\tz\eta}|2^k\rho(x)|^{n}}\dint_{|x-y|<2^k\rho(x)}|b(x)-b(y)||f(y)|dy\\
&\le CM_{V,\eta}^bf(x),
\end{array}\eqno(4.7)$$
if $N>M>n+\tz\eta$.

Thus, by (4.6) and (4.7), and using Lemmas 3.3, 3.4, 3.5, Theorem 2.2 and Proposition 2.2, we can obtain the desired result.                            \hfill$\Box$

We remark that in fact all results in this section also hold for $BMO_{\tz_1}(\rho)$ and $A_p^{\rho,\tz_2}$ if $\tz_1\not=\tz_2$.

 LMAM, School of Mathematical   Sciences

 Peking University

 Beijing, 100871

 P. R. China

\bigskip

 E-mail address:  tanglin@math.pku.edu.cn

\end{document}